\documentclass[10pt]{article}
\usepackage{amsmath}
\usepackage{amssymb}
\usepackage{amsthm}
\usepackage{mathrsfs}

\begin{document}

\title{Boundedness of fractional integral operators with rough kernels on weighted Morrey spaces}
\author{Hua Wang \footnote{E-mail address: wanghua@pku.edu.cn.}\\
\footnotesize{Department of Mathematics, Zhejiang University, Hangzhou 310027, China}}
\date{}
\maketitle

\begin{abstract}
Let $M_{\Omega,\alpha}$ and $T_{\Omega,\alpha}$ be the fractional maximal and integral operators with rough kernels, where $0<\alpha<n$. In this paper, we shall study the continuity properties of $M_{\Omega,\alpha}$ and $T_{\Omega,\alpha}$ on the weighted Morrey spaces $L^{p,\kappa}(w)$. The boundedness of their commutators with BMO functions is also obtained.
\end{abstract}

\section{Introduction}

Let $\Omega\in L^s(S^{n-1})$ be homogeneous of degree zero on $\mathbb R^n$, where $S^{n-1}$ denotes the unit sphere of $\mathbb R^n$($n\ge2$) equipped with the normalized Lebesgue measure $d\sigma$ and $s>1$. For any $0<\alpha<n$, then the fractional integral operator with rough kernel $T_{\Omega,\alpha}$ is defined by
\begin{equation*}
T_{\Omega,\alpha}f(x)=\int_{\mathbb R^n}\frac{\Omega(y')}{|y|^{n-\alpha}}f(x-y)\,dy
\end{equation*}
and a related fractional maximal operator $M_{\Omega,\alpha}$ is defined by
\begin{equation*}
M_{\Omega,\alpha}f(x)=\sup_{r>0}\frac{1}{r^{n-\alpha}}\int_{|y|\le r}\big|\Omega(y')f(x-y)\big|\,dy,
\end{equation*}
where $y'=y/{|y|}$ for any $y\neq0$. In 1971, Muckenhoupt and Wheeden \cite{muckenhoupt3} studied the weighted norm inequalities for $T_{\Omega,\alpha}$ with the weight $w(x)=|x|^\beta$. The weak type estimates with power weights for $M_{\Omega,\alpha}$ and $T_{\Omega,\alpha}$ was obtained by Ding in \cite{ding3}. Later, Ding and Lu \cite{ding1} considered the weighted norm inequalities for $M_{\Omega,\alpha}$ and $T_{\Omega,\alpha}$ with more general weights. More precisely, they proved

\newtheorem*{thma}{Theorem A}
\begin{thma}[\cite{ding1}]
Let $0<\alpha<n$, $1\le s'<p<n/{\alpha}$ and $1/q=1/p-{\alpha}/n$. If $\Omega\in L^s(S^{n-1})$ and $w^{s'}\in A(p/{s'},q/{s'})$, then the operators $M_{\Omega,\alpha}$ and $T_{\Omega,\alpha}$ are all bounded from $L^p(w^p)$ to $L^q(w^q)$.
\end{thma}

Let $b$ be a locally integrable function on $\mathbb R^n$, then for $0<\alpha<n$, we shall define the commutators generated by fractional maximal and integral operators with rough kernels and $b$ as follows.
\begin{equation*}
[b,M_{\Omega,\alpha}](f)(x)=\sup_{r>0}\frac{1}{r^{n-\alpha}}\int_{|y-x|\le r}|b(x)-b(y)||\Omega(x-y)f(y)|\,dy,
\end{equation*}
\begin{equation*}
\begin{split}
[b,T_{\Omega,\alpha}](f)(x)&=b(x)T_{\Omega,\alpha}f(x)-T_{\Omega,\alpha}(bf)(x)\\
&=\int_{\mathbb R^n}\frac{\Omega(x-y)}{|x-y|^{n-\alpha}}[b(x)-b(y)]f(y)\,dy.
\end{split}
\end{equation*}

In 1993, by using the Rubio de Francia extrapolation theorem, Segovia and Torrea \cite{segovia2} obtained the weighted boundedness of commutator $[b,T_{\Omega,\alpha}]$, where $b\in BMO(\mathbb R^n)$ and $\Omega$ satisfies some Dini smoothness condition (see also \cite{segovia1}). In 1999, Ding and Lu \cite{ding2} improved this result by removing the smoothness condition imposed on $\Omega$. More specifically, they showed (see also \cite{lu}).

\newtheorem*{thmb}{Theorem B}
\begin{thmb}[\cite{ding2}]
Let $0<\alpha<n$, $1\le s'<p<n/{\alpha}$ and $1/q=1/p-{\alpha}/n$. Assume that $\Omega\in L^s(S^{n-1})$, $w^{s'}\in A(p/{s'},q/{s'})$ and $b\in BMO(\mathbb R^n)$, then the commutator $[b,T_{\Omega,\alpha}]$ is bounded from $L^p(w^p)$ to $L^q(w^q)$.
\end{thmb}

The classical Morrey spaces $\mathcal L^{p,\lambda}$ were first introduced by Morrey in \cite{morrey} to study the local behavior of solutions to second order elliptic partial differential equations. For the boundedness of the Hardy-Littlewood maximal operator, the fractional integral operator and the Calder\'on-Zygmund singular integral operator on these spaces, we refer the readers to \cite{adams,chiarenza,peetre}. For the properties and applications of classical Morrey spaces, see \cite{fan,fazio1,fazio2} and references therein.

In 2009, Komori and Shirai \cite{komori} first defined the weighted Morrey spaces $L^{p,\kappa}(w)$ which could be viewed as an extension of weighted Lebesgue spaces, and studied the boundedness of the above classical operators on these weighted spaces. Recently, in \cite{wang1} and \cite{wang2}, we have established the continuity properties of some other operators on the weighted Morrey spaces $L^{p,\kappa}(w)$.

The purpose of this paper is to discuss the boundedness properties of $M_{\Omega,\alpha}$ and $T_{\Omega,\alpha}$ on the weighted Morrey spaces. Here, and in what follows we shall use the notation $s'=s/{(s-1)}$ when $1<s<\infty$ and $s'=1$ when $s=\infty$. Our main results in the paper are formulated as follows.

\newtheorem{theorem}{Theorem}[section]

\begin{theorem}
Suppose that $\Omega\in L^s(S^{n-1})$ with $1<s\le\infty$. If $0<\alpha<n$, $1\le s'<p<n/{\alpha}$, $1/q=1/p-{\alpha}/n$, $0<\kappa<p/q$ and $w^{s'}\in A(p/{s'},q/{s'})$, then the fractional maximal operator $M_{\Omega,\alpha}$ is bounded from $L^{p,\kappa}(w^p,w^q)$ to $L^{q,{\kappa q}/p}(w^q)$.
\end{theorem}

\begin{theorem}
Suppose that $\Omega\in L^s(S^{n-1})$ with $1<s\le\infty$. If $0<\alpha<n$, $1\le s'<p<n/{\alpha}$, $1/q=1/p-{\alpha}/n$, $0<\kappa<p/q$ and $w^{s'}\in A(p/{s'},q/{s'})$, then the fractional integral operator $T_{\Omega,\alpha}$ is bounded from $L^{p,\kappa}(w^p,w^q)$ to $L^{q,{\kappa q}/p}(w^q)$.
\end{theorem}

\begin{theorem}
Suppose that $\Omega\in L^s(S^{n-1})$ with $1<s\le\infty$ and $b\in BMO(\mathbb R^n)$. If $0<\alpha<n$, $1\le s'<p<n/{\alpha}$, $1/q=1/p-{\alpha}/n$, $0<\kappa<p/q$ and $w^{s'}\in A(p/{s'},q/{s'})$, then the commutator $[b,T_{\Omega,\alpha}]$ is bounded from $L^{p,\kappa}(w^p,w^q)$ to $L^{q,{\kappa q}/p}(w^q)$.
\end{theorem}

\section{Notations and definitions}

Let us first recall some standard definitions and notations. The classical $A_p$ weight theory was first introduced by Muckenhoupt in the study of weighted $L^p$ boundedness of Hardy-Littlewood maximal functions in \cite{muckenhoupt1}. A weight $w$ is a nonnegative, locally integrable function on $\mathbb R^n$, $B=B(x_0,r_B)$ denotes the ball
with the center $x_0$ and radius $r_B$. Given a ball $B$ and $\lambda>0$, $\lambda B$ denotes the ball with the same center as $B$ whose radius is $\lambda$ times that of $B$. For a given weight function $w$, we also denote
the Lebesgue measure of $B$ by $|B|$ and the weighted measure of $B$ by $w(B)$, where $w(B)=\int_B w(x)\,dx$. We say that $w\in A_p$, $1<p<\infty$, if
\begin{equation*}
\left(\frac1{|B|}\int_B w(x)\,dx\right)\left(\frac1{|B|}\int_B w(x)^{-1/{(p-1)}}\,dx\right)^{p-1}\le C \quad\mbox{for every ball}\; B\subseteq \mathbb
R^n,
\end{equation*}
where $C$ is a positive constant which is independent of $B$.

For the case $p=1$, $w\in A_1$, if
\begin{equation*}
\frac1{|B|}\int_B w(x)\,dx\le C\cdot\underset{x\in B}{\mbox{ess\,inf}}\,w(x)\quad\mbox{for every ball}\;B\subseteq\mathbb R^n.
\end{equation*}

For the case $p=\infty$, $w\in A_\infty$ if it satisfies the $A_p$ condition for some $1<p<\infty$.

We also need another weight class $A(p,q)$ introduced by Muckenhoupt and Wheeden in \cite{muckenhoupt2}. A weight function $w$ belongs to $A(p,q)$ for $1<p<q<\infty$ if there exists a constant $C>0$ such that
\begin{equation*}
\left(\frac{1}{|B|}\int_B w(x)^q\,dx\right)^{1/q}\left(\frac{1}{|B|}\int_B w(x)^{-p'}\,dx\right)^{1/{p'}}\le C \quad\mbox{for every ball}\; B\subseteq \mathbb R^n.
\end{equation*}

A weight function $w$ is said to belong to the reverse H\"{o}lder class $RH_r$ if there exist two constants $r>1$ and $C>0$ such that the following reverse H\"{o}lder inequality holds
\begin{equation*}
\left(\frac{1}{|B|}\int_B w(x)^r\,dx\right)^{1/r}\le C\left(\frac{1}{|B|}\int_B w(x)\,dx\right)\quad\mbox{for every ball}\; B\subseteq \mathbb R^n.
\end{equation*}

We state the following results that we will use frequently in the sequel.

\newtheorem{lemma}[theorem]{Lemma}
\begin{lemma}[\cite{garcia}]
Let $w\in A_p$ with $p\ge1$. Then, for any ball $B$, there exists an absolute constant $C>0$ such that
\begin{equation*}
w(2B)\le C\,w(B).
\end{equation*}
In general, for any $\lambda>1$, we have
\begin{equation*}
w(\lambda B)\le C\cdot\lambda^{np}w(B),
\end{equation*}
where $C$ does not depend on $B$ nor on $\lambda$.
\end{lemma}

\begin{lemma}[\cite{gundy}]
Let $w\in RH_r$ with $r>1$. Then there exists a constant $C>0$ such that
\begin{equation*}
\frac{w(E)}{w(B)}\le C\left(\frac{|E|}{|B|}\right)^{(r-1)/r}
\end{equation*}
for any measurable subset $E$ of a ball $B$.
\end{lemma}

Next we shall introduce the Hardy-Littlewood maximal operator, its variant and BMO spaces.
The Hardy-Littlewood maximal operator $M$ is defined by
\begin{equation*}
M(f)(x)=\sup_{x\in B}\frac{1}{|B|}\int_B|f(y)|\,dy, \end{equation*}
where the supremum is taken over all balls $B$ containing $x$. For $0<\alpha<n$, $s\ge1$, we define the fractional maximal operator $M_{\alpha,s}$ by
\begin{equation*}
M_{\alpha,s}(f)(x)=\sup_{x\in B}\bigg(\frac{1}{|B|^{1-\frac{\alpha s}{n}}}\int_B|f(y)|^s\,dy\bigg)^{1/s}.
\end{equation*}
Moreover, we denote simply by $M_\alpha$ when $s=1$.

A locally integrable function $b$ is said to be in $BMO(\mathbb R^n)$ if
\begin{equation*}
\|b\|_*=\sup_{B}\frac{1}{|B|}\int_B|b(x)-b_B|\,dx<\infty,
\end{equation*}
where $b_B$ stands for the average of $b$ on $B$, i.e. $b_B=\frac{1}{|B|}\int_B b(y)\,dy$ and the supremum is taken over all balls $B$ in $\mathbb R^n$.

\newtheorem*{thmc}{Theorem C}
\begin{thmc}[\cite{duoand,john}]
Assume that $b\in BMO(\mathbb R^n)$. Then for any $1\le p<\infty$, we have
\begin{equation*}
\sup_B\bigg(\frac{1}{|B|}\int_B\big|b(x)-b_B\big|^p\,dx\bigg)^{1/p}\le C\|b\|_*.
\end{equation*}
\end{thmc}

We are going to conclude this section by defining the weighted Morrey space
and giving the known result relevant to this paper. For further details, we refer the readers to \cite{komori}.
\newtheorem{defn}[theorem]{Definition}

\begin{defn}[\cite{komori}]
Let $1\le p<\infty$, $0<\kappa<1$ and $w$ be a weight function. Then the weighted Morrey space is defined by
\begin{equation*}
L^{p,\kappa}(w)=\big\{f\in L^p_{loc}(w):\|f\|_{L^{p,\kappa}(w)}<\infty\big\},
\end{equation*}
where
\begin{equation*}
\|f\|_{L^{p,\kappa}(w)}=\sup_B\left(\frac{1}{w(B)^\kappa}\int_B|f(x)|^pw(x)\,dx\right)^{1/p}
\end{equation*}
and the supremum is taken over all balls $B$ in $\mathbb R^n$.
\end{defn}
In order to deal with the fractional order case, we need to consider the weighted Morrey space with two weights.

\begin{defn}[\cite{komori}]
Let $1\le p<\infty$ and $0<\kappa<1$. Then for two weights $u$ and $v$, the weighted Morrey space is defined by
\begin{equation*}
L^{p,\kappa}(u,v)=\big\{f\in L^p_{loc}(u):\|f\|_{L^{p,\kappa}(u,v)}<\infty\big\},
\end{equation*}
where
\begin{equation*}
\|f\|_{L^{p,\kappa}(u,v)}=\sup_{B}\left(\frac{1}{v(B)^{\kappa}}\int_B|f(x)|^pu(x)\,dx\right)^{1/p}.
\end{equation*}
\end{defn}

\newtheorem*{thmd}{Theorem D}
\begin{thmd}
If $0<\alpha<n$, $1<p<n/{\alpha}$, $1/q=1/p-\alpha/n$, $0<\kappa<p/q$ and $w\in A(p,q)$, then the fractional maximal operator $M_\alpha$ is bounded from $L^{p,\kappa}(w^p,w^q)$ to $L^{q,{\kappa q}/p}(w^q)$.
\end{thmd}

Throughout this article, we will use $C$ to denote a positive constant, which is independent of the main parameters and not necessarily the same at each occurrence. By $A\sim B$, we mean that there exists a constant $C>1$ such that $\frac1C\le\frac AB\le C$.

\section{Proof of Theorem 1.1}

\begin{proof}[Proof of Theorem 1.1]
For $\Omega\in L^s(S^{n-1})$, we set
\begin{equation*}
\|\Omega\|_{L^s(S^{n-1})}=\bigg(\int_{S^{n-1}}\big|\Omega(y')\big|^s\,d\sigma(y')\bigg)^{1/s}.
\end{equation*}
From H\"older's inequality, it follows that
\begin{equation*}
\begin{split}
M_{\Omega,\alpha}f(x)&\le\sup_{r>0}\frac{1}{r^{n-\alpha}}\bigg(\int_{|y|\le r}\big|\Omega(y')\big|^s\,dy\bigg)^{1/s}\bigg(\int_{|y|\le r}|f(x-y)|^{s'}\,dy\bigg)^{1/{s'}}\\
&\le C\|\Omega\|_{L^s(S^{n-1})}\cdot\sup_{r>0}\bigg(\frac{1}{r^{n-\alpha s'}}\int_{|y|\le r}|f(x-y)|^{s'}\,dy\bigg)^{1/{s'}}\\
&\le C\|\Omega\|_{L^s(S^{n-1})}\cdot\sup_{r>0}\Bigg(\frac{1}{|B(x,r)|^{1-\frac{\alpha s'}{n}}}\int_{B(x,r)}|f(y)|^{s'}\,dy\Bigg)^{1/{s'}}\\
&= C\|\Omega\|_{L^s(S^{n-1})}M_{\alpha,s'}(f)(x).
\end{split}
\end{equation*}
If we let $p_1=p/{s'}$, $q_1=q/{s'}$ and $\nu=w^{s'}$, then for $0<\alpha<n$, $1\le s'<n/{\alpha}$, we have $1/{q_1}=1/{p_1}-{(\alpha s')}/n$ and $0<\kappa<{p_1}/{q_1}$. Also observe that
\begin{equation*}
M_{\alpha,s'}(f)=M_{\alpha s'}(|f|^{s'})^{1/{s'}}.
\end{equation*}
Hence, by Theorem D, we obtain
\begin{equation*}
\begin{split}
\big\|M_{\alpha,s'}(f)\big\|_{L^{q,{\kappa q}/p}(w^q)}&=\big\|M_{\alpha s'}(|f|^{s'})\big\|^{1/{s'}}_{L^{q_1,{\kappa q_1}/{p_1}}(\nu^{q_1})}\\
&\le C\big\||f|^{s'}\big\|^{1/{s'}}_{L^{p_1,\kappa}(\nu^{p_1},\nu^{q_1})}\\
&\le C\|f\|_{L^{p,\kappa}(w^p,w^q)}.
\end{split}
\end{equation*}
This finishes the proof of Theorem 1.1.
\end{proof}

\section{Proof of Theorem 1.2}

\begin{proof}[Proof of Theorem 1.2]
Fix a ball $B=B(x_0,r_B)\subseteq\mathbb R^n$ and decompose $f=f_1+f_2$,
where $f_1=f\chi_{_{2B}}$ and $\chi_{_{2B}}$ denotes the characteristic function of $2B$.
Since $T_{\Omega,\alpha}$ is a linear operator, then we can write
\begin{equation*}
\begin{split}
&\frac{1}{w^q(B)^{\kappa/p}}\bigg(\int_B|T_{\Omega,\alpha}f(x)|^qw(x)^q\,dx\bigg)^{1/q}\\
\le&\,\frac{1}{w^q(B)^{\kappa/p}}\bigg(\int_B|T_{\Omega,\alpha}f_1(x)|^qw(x)^q\,dx\bigg)^{1/q}
+\frac{1}{w^q(B)^{\kappa/p}}\bigg(\int_B|T_{\Omega,\alpha}f_2(x)|^qw(x)^q\,dx\bigg)^{1/q}\\
=&\,I_1+I_2.
\end{split}
\end{equation*}
As in the proof of Theorem 1.1, we also set $p_1=p/{s'}$, $q_1=q/{s'}$ and $\nu=w^{s'}$.
Since $\nu\in A(p_1,q_1)$, then we get $\nu^{q_1}=w^q\in A_{1+{q_1}/{p'_1}}$ (see \cite{muckenhoupt2}).
Hence, by Theorem A and Lemma 2.1, we have
\begin{equation*}
\begin{split}
I_1&\le C\cdot\frac{1}{w^q(B)^{\kappa/p}}\bigg(\int_{2B}|f(x)|^pw(x)^p\,dx\bigg)^{1/p}\\
&\le C\|f\|_{L^{p,\kappa}(w^p,w^q)}\cdot\frac{w^q(2B)^{\kappa/p}}{w^q(B)^{\kappa/p}}\\
&\le C\|f\|_{L^{p,\kappa}(w^p,w^q)}.
\end{split}
\end{equation*}

We now turn to deal with the term $I_2$. An application of H\"older's inequality gives us that
\begin{align}
\big|T_{\Omega,\alpha}(f_2)(x)\big|&\le\int_{(2B)^c}\frac{|\Omega(x-y)|}{|x-y|^{n-\alpha}}|f(y)|\,dy\\
&\le\sum_{k=1}^\infty\bigg(\int_{2^{k+1}B\backslash 2^{k}B}|\Omega(x-y)|^s\,dy\bigg)^{1/s}\bigg(\int_{2^{k+1}B\backslash 2^{k}B}\frac{|f(y)|^{s'}}{|x-y|^{(n-\alpha)s'}}\,dy\bigg)^{1/{s'}}\notag.
\end{align}
When $x\in B$ and $y\in2^{k+1}B\backslash 2^{k}B$, then we can easily see that
$2^{k-1}r_B\le|y-x|<2^{k+2}r_B$. Thus, by a simple computation, we deduce
\begin{equation}
\bigg(\int_{2^{k+1}B\backslash 2^{k}B}|\Omega(x-y)|^s\,dy\bigg)^{1/s}\le C\|\Omega\|_{L^s(S^{n-1})}\big|2^{k+1}B\big|^{1/s}.
\end{equation}
We also note that if $x\in B$, $y\in(2B)^c$, then $|y-x|\sim|y-x_0|$. Consequently
\begin{equation}
\bigg(\int_{2^{k+1}B\backslash 2^{k}B}\frac{|f(y)|^{s'}}{|x-y|^{(n-\alpha)s'}}\,dy\bigg)^{1/{s'}}\le C\cdot\frac{1}{|2^{k+1}B|^{1-\alpha/n}}\bigg(\int_{2^{k+1}B}|f(y)|^{s'}\,dy\bigg)^{1/{s'}}.
\end{equation}
Substituting the above two inequalities (2) and (3) into (1), we obtain
\begin{equation*}
\big|T_{\Omega,\alpha}(f_2)(x)\big|\le C\|\Omega\|_{L^s(S^{n-1})}\sum_{k=1}^\infty
\frac{1}{|2^{k+1}B|^{1-\alpha/n-1/s}}
\bigg(\int_{2^{k+1}B}|f(y)|^{s'}\,dy\bigg)^{1/{s'}}.
\end{equation*}
By using H\"{o}lder's inequality and the definition of $\nu\in A(p_1,q_1)$, we can get
\begin{align}
\bigg(\int_{2^{k+1}B}|f(y)|^{s'}\,dy\bigg)^{1/{s'}}&\le
\bigg(\int_{2^{k+1}B}|f(y)|^{p_1s'}\nu(y)^{p_1}\,dy\bigg)^{1/{(p_1s')}}
\bigg(\int_{2^{k+1}B}\nu(y)^{-p'_1}\,dy\bigg)^{1/{(p'_1s')}}\notag\\
&\le C\bigg(\int_{2^{k+1}B}|f(y)|^{p}w(y)^{p}\,dy\bigg)^{1/{p}}
\bigg(\frac{|2^{k+1}B|^{1-1/{p_1}+1/{q_1}}}{\nu^{q_1}(2^{k+1}B)^{1/{q_1}}}\bigg)^{1/{s'}}\notag\\
&\le C\|f\|_{L^{p,\kappa}(w^p,w^q)}w^q\big(2^{k+1}B\big)^{\kappa/p}
\cdot\frac{|2^{k+1}B|^{1/{s'}-1/p+1/q}}{w^q(2^{k+1}B)^{1/q}}\notag\\
&= C\|f\|_{L^{p,\kappa}(w^p,w^q)}\big|2^{k+1}B\big|^{1-1/s-\alpha/n}\cdot w^q\big(2^{k+1}B\big)^{\kappa/p-1/q}.
\end{align}
So we have
\begin{equation*}
\big|T_{\Omega,\alpha}(f_2)(x)\big|\le C \|f\|_{L^{p,\kappa}(w^p,w^q)}
\sum_{k=1}^\infty w^q\big(2^{k+1}B\big)^{\kappa/p-1/q},
\end{equation*}
which implies
\begin{equation*}
I_2\le C\|f\|_{L^{p,\kappa}(w^p,w^q)}\sum_{k=1}^\infty\frac{w^q(B)^{1/q-\kappa/p}}{w^q(2^{k+1}B)^{1/q-\kappa/p}}.
\end{equation*}
Observe that $w^q=\nu^{q_1}\in A_{1+{q_1}/{p'_1}}$, then we know that there exists $r>1$ such that $w^q\in RH_r$.
Thus, it follows directly from Lemma 2.2 that
\begin{equation}
\frac{w^q(B)}{w^q(2^{k+1}B)}\le C\left(\frac{|B|}{|2^{k+1}B|}\right)^{{(r-1)}/r}.
\end{equation}
Therefore
\begin{equation*}
\begin{split}
I_2&\le C\|f\|_{L^{p,\kappa}(w^p,w^q)}\sum_{k=1}^\infty\left(\frac{1}{2^{kn}}\right)^{(1-1/r)(1/q-\kappa/p)}\\
&\le C\|f\|_{L^{p,\kappa}(w^p,w^q)},
\end{split}
\end{equation*}
where the last series is convergent since $r>1$ and $0<\kappa<p/q$.
Combining the above estimates for $I_1$ and $I_2$
and taking the supremum over all balls $B\subseteq\mathbb R^n$, we complete the proof of Theorem 1.2.
\end{proof}

\section{Proof of Theorem 1.3}

\begin{proof}[Proof of Theorem 1.3]
Fix a ball $B=B(x_0,r_B)\subseteq\mathbb R^n$. Let $f=f_1+f_2$, where $f_1=f\chi_{_{2B}}$.
Since $[b,T_{\Omega,\alpha}]$ is a linear operator, then we have
\begin{equation*}
\begin{split}
&\frac{1}{w^q(B)^{\kappa/p}}\bigg(\int_B\big|[b,T_{\Omega,\alpha}]f(x)\big|^qw(x)^q\,dx\bigg)^{1/q}\\
\le&\,\frac{1}{w^q(B)^{\kappa/p}}\bigg(\int_B\big|[b,T_{\Omega,\alpha}]f_1(x)\big|^qw(x)^q\,dx\bigg)^{1/q}
+\frac{1}{w^q(B)^{\kappa/p}}\bigg(\int_B\big|[b,T_{\Omega,\alpha}]f_2(x)\big|^qw(x)^q\,dx\bigg)^{1/q}\\
=&\,J_1+J_2.
\end{split}
\end{equation*}
As before, we set $p_1=p/{s'}$, $q_1=q/{s'}$ and $\nu=w^{s'}$, then $\nu^{q_1}=w^q\in A_{1+{q_1}/{p'_1}}$.
Theorem B and Lemma 2.1 imply
\begin{align} J_1&\le C\|b\|_*\cdot\frac{1}{w^q(B)^{\kappa/p}}\bigg(\int_{2B}|f(x)|^pw(x)^p\,dx\bigg)^{1/p}\notag\\
&\le C\|b\|_*\|f\|_{L^{p,\kappa}(w^p,w^q)}\cdot\frac{w^q(2B)^{\kappa/p}}{w^q(B)^{\kappa/p}}\notag\\
&\le C\|b\|_*\|f\|_{L^{p,\kappa}(w^p,w^q)}. \end{align}

In order to estimate the term $J_2$, for any $x\in B$, we first write
\begin{equation*}
\begin{split}
\big|\big[b,T_{\Omega,\alpha}\big]f_2(x)\big|=&\,\bigg|\int_{(2B)^c}\frac{\Omega(x-y)}{|x-y|^{n-\alpha}}
\big[b(x)-b(y)\big]f(y)\,dy\bigg|\\
\le&\,\big|b(x)-b_B\big|\cdot\int_{(2B)^c}\frac{|\Omega(x-y)|}{|x-y|^{n-\alpha}}|f(y)|\,dy\\
&+\int_{(2B)^c}\frac{|\Omega(x-y)|}{|x-y|^{n-\alpha}}|b(y)-b_B||f(y)|\,dy\\
=&\,\mbox{\upshape I+II}.
\end{split}
\end{equation*}
For the term I, it follows from the previous estimates (2) and (4) that
\begin{equation*}
\mbox{\upshape I}\le C\|f\|_{L^{p,\kappa}(w^p,w^q)}|b(x)-b_B|\cdot\sum_{k=1}^\infty
\frac{1}{w^q(2^{k+1}B)^{1/q-\kappa/p}}.
\end{equation*}
Hence
\begin{equation*}
\begin{split}
&\frac{1}{w^q(B)^{\kappa/p}}\bigg(\int_B\mbox{\upshape I}^q\,w(x)^q\,dx\bigg)^{1/q}\\
\le &\,
C\|f\|_{L^{p,\kappa}(w^p,w^q)}\frac{1}{w^q(B)^{\kappa/p}}
\cdot\sum_{k=1}^\infty\frac{1}{w^q(2^{k+1}B)^{1/q-\kappa/p}}\cdot\bigg(\int_B|b(x)-b_B|^qw(x)^q\,dx\bigg)^{1/q}\\
=&\,C\|f\|_{L^{p,\kappa}(w^p,w^q)}\sum_{k=1}^\infty\frac{w^q(B)^{1/q-\kappa/p}}{w^q(2^{k+1}B)^{1/q-\kappa/p}}
\cdot\bigg(\frac{1}{w^q(B)}\int_B|b(x)-b_B|^qw(x)^q\,dx\bigg)^{1/q}.
\end{split}
\end{equation*}
We now claim that for any $1<q<\infty$ and $\mu\in A_\infty$, the following inequality holds
\begin{equation}
\bigg(\frac{1}{\mu(B)}\int_B|b(x)-b_B|^q\mu(x)\,dx\bigg)^{1/q}\le C\|b\|_*.
\end{equation}
In fact, since $\mu\in A_\infty$, then there must exist $r>1$ such that $\mu\in RH_r$.
Thus, by H\"older's inequality and Theorem C, we obtain
\begin{equation*}
\begin{split}
\bigg(\frac{1}{\mu(B)}\int_B|b(x)-b_B|^q\mu(x)\,dx\bigg)^{1/q}&\le\frac{1}{\mu(B)^{1/q}}
\bigg(\int_B|b(x)-b_B|^{qr'}\,dx\bigg)^{1/{(qr')}}\bigg(\int_B \mu(x)^r\,dx\bigg)^{1/{(qr)}}\\
&\le C\bigg(\frac{1}{|B|}\int_B|b(x)-b_B|^{qr'}\,dx\bigg)^{1/{(qr')}}\\
&\le C\|b\|_*,
\end{split}
\end{equation*}
which is our desired result. Note that $w^q\in A_{1+{q_1}/{p'_1}}\subset A_\infty$.
In addition, we have $w^q\in RH_r$ with
$r>1$. Hence, by the inequalities (5) and (7), we get
\begin{align}
\frac{1}{w^q(B)^{\kappa/p}}\bigg(\int_B \mbox{\upshape I}^q\,w(x)^q\,dx\bigg)^{1/q}
&\le
C\|b\|_*\|f\|_{L^{p,\kappa}(w^p,w^q)}\sum_{k=1}^\infty\left(\frac{1}{2^{kn}}\right)^{(1-1/r)(1/q-\kappa/p)}\notag\\
&\le C\|b\|_*\|f\|_{L^{p,\kappa}(w^p,w^q)}.
\end{align}
On the other hand
\begin{equation*}
\begin{split}
\mbox{\upshape II}\le&\,\sum_{k=1}^\infty\int_{2^{k+1}B\backslash 2^kB}\frac{|\Omega(x-y)|}{|x-y|^{n-\alpha}}
|b(y)-b_B||f(y)|\,dy\\
\le&\,\sum_{k=1}^\infty\int_{2^{k+1}B\backslash 2^kB}\frac{|\Omega(x-y)|}{|x-y|^{n-\alpha}}
\big|b(y)-b_{2^{k+1}B}\big||f(y)|\,dy\\
&+\sum_{k=1}^\infty\int_{2^{k+1}B\backslash 2^kB}\frac{|\Omega(x-y)|}{|x-y|^{n-\alpha}}
\big|b_{2^{k+1}B}-b_B\big||f(y)|\,dy\\
=&\,\mbox{\upshape III+IV}.
\end{split}
\end{equation*}
To estimate III and IV, we observe that when $x\in B$, $y\in (2B)^c$, then $|y-x|\sim|y-x_0|$.
Thus, it follows from H\"older's inequality and (2) that
\begin{equation*}
\begin{split}
\mbox{\upshape III}\le C \sum_{k=1}^\infty\frac{1}{|2^{k+1}B|^{1-\alpha/n}}
\cdot\int_{2^{k+1}B\backslash 2^kB}\big|\Omega(x-y)\big|\big|b(y)-b_{2^{k+1}B}\big||f(y)|\,dy\\
\le C \sum_{k=1}^\infty\frac{1}{|2^{k+1}B|^{1-\alpha/n-1/s}}
\cdot\bigg(\int_{2^{k+1}B}\big|b(y)-b_{2^{k+1}B}\big|^{s'}|f(y)|^{s'}\,dy\bigg)^{1/{s'}}.
\end{split}
\end{equation*}
An application of H\"older's inequality yields
\begin{equation*}
\begin{split}
&\bigg(\int_{2^{k+1}B}\big|b(y)-b_{2^{k+1}B}\big|^{s'}|f(y)|^{s'}\,dy\bigg)^{1/{s'}} \\
\le&\,
\bigg(\int_{2^{k+1}B}|f(y)|^{p_1s'}\nu(y)^{p_1}\,dy\bigg)^{1/{(p_1s')}}
\bigg(\int_{2^{k+1}B}\big|b(y)-b_{2^{k+1}B}\big|^{p'_1s'}\nu(y)^{-p'_1}\,dy\bigg)^{1/{(p'_1s')}}\\
\le&\, \bigg(\int_{2^{k+1}B}|f(y)|^{p}w(y)^{p}\,dy\bigg)^{1/{p}}
\bigg(\int_{2^{k+1}B}\big|b(y)-b_{2^{k+1}B}\big|^{p'_1s'}\nu(y)^{-p'_1}\,dy\bigg)^{1/{(p'_1s')}}.
\end{split}
\end{equation*}
Since $\nu\in A(p_1,q_1)$, then we know that $\nu^{-p'_1}\in A_{1+{p'_1}/{q_1}}\subset A_\infty$(see \cite{muckenhoupt2}). Hence, by using the inequality (7) and the fact that $\nu\in A(p_1,q_1)$, we obtain
\begin{align*}
\bigg(\int_{2^{k+1}B}\big|b(y)-b_{2^{k+1}B}\big|^{p'_1s'}\nu(y)^{-p'_1}\,dy\bigg)^{1/{(p'_1s')}}
&\le C\|b\|_*\cdot\nu^{-p'_1}\big(2^{k+1}B\big)^{1/{(p'_1s')}} \notag\\
&\le C\|b\|_*\cdot\left(\frac{|2^{k+1}B|^{1/{q_1}+1/{p'_1}}}{\nu^{q_1}(2^{k+1}B)^{1/{q_1}}}\right)^{1/{s'}}\notag\\
\end{align*}
\begin{align}
&= C\|b\|_*\cdot\frac{|2^{k+1}B|^{1/{s'}-1/p+1/q}}{w^q(2^{k+1}B)^{1/q}}.
\end{align}
Consequently, by the above inequality (9), we deduce
\begin{equation*}
\mbox{\upshape III}\le C\|b\|_*\|f\|_{L^{p,\kappa}(w^p,w^q)}
\sum_{k=1}^\infty \frac{1}{w^q(2^{k+1}B)^{1/q-\kappa/p}},
\end{equation*}
which implies
\begin{align}
\frac{1}{w^q(B)^{\kappa/p}}\bigg(\int_B \mbox{\upshape III}^q\,w(x)^q\,dx\bigg)^{1/q}
&\le C\|b\|_*\|f\|_{L^{p,\kappa}(w^p,w^q)}
\sum_{k=1}^\infty\frac{w^q(B)^{1/q-\kappa/p}}{w^q(2^{k+1}B)^{1/q-\kappa/p}}\notag\\
&\le C\|b\|_*\|f\|_{L^{p,\kappa}(w^p,w^q)}.
\end{align}
Since $b\in BMO(\mathbb R^n)$, then a direct calculation shows that
\begin{equation*}
\big|b_{2^{k+1}B}-b_B\big|\le C\cdot k\|b\|_*.
\end{equation*}
Moreover, by H\"older's inequality, the estimates (2) and (4), we can get
\begin{equation*}
\begin{split}
\mbox{\upshape IV}&\le C\|b\|_*\sum_{k=1}^\infty k\cdot\frac{1}{|2^{k+1}B|^{1-\alpha/n}}
\int_{2^{k+1}B\backslash 2^kB}|\Omega(x-y)||f(y)|\,dy\\
&\le C\|b\|_*\|f\|_{L^{p,\kappa}(w^p,w^q)}\sum_{k=1}^\infty k\cdot\frac{1}{w^q(2^{k+1}B)^{1/q-\kappa/p}}.
\end{split}
\end{equation*}
Therefore
\begin{align}
\frac{1}{w^q(B)^{\kappa/p}}\bigg(\int_B \mbox{\upshape IV}^q\,w(x)^q\,dx\bigg)^{1/q}
&\le C\|b\|_*\|f\|_{L^{p,\kappa}(w^p,w^q)}
\sum_{k=1}^\infty k\cdot\frac{w^q(B)^{1/q-\kappa/p}}{w^q(2^{k+1}B)^{1/q-\kappa/p}}\notag\\
&\le C\|b\|_*\|f\|_{L^{p,\kappa}(w^p,w^q)}\sum_{k=1}^\infty \frac{k}{2^{kn\delta}}\notag\\
&\le C\|b\|_*\|f\|_{L^{p,\kappa}(w^p,w^q)},
\end{align}
where $w^q\in RH_r$ and $\delta=(1-1/r)(1/q-\kappa/p)$. Summarizing the estimates (10)
and (11) derived above, we thus obtain
\begin{equation}
\frac{1}{w^q(B)^{\kappa/p}}\bigg(\int_B \mbox{\upshape II}^q\,w(x)^q\,dx\bigg)^{1/q}
\le C\|b\|_*\|f\|_{L^{p,\kappa}(w^p,w^q)}.
\end{equation}
Combining the inequalities (6) and (8) with the above inequality (12) and
taking the supremum over all balls $B\subseteq\mathbb R^n$, we conclude the proof of Theorem 1.3.
\end{proof}

It should be pointed out that $[b,M_{\Omega,\alpha}](f)$ can be controlled pointwise
by $[b,T_{|\Omega|,\alpha}](|f|)$ for any $f(x)$. In fact, for any $0<\alpha<n$, $x\in\mathbb R^n$ and $r>0$, we have
\begin{equation*}
\begin{split}
[b,T_{|\Omega|,\alpha}](|f|)(x)&\ge\int_{|y-x|\le r}\frac{|\Omega(x-y)|}{|x-y|^{n-\alpha}}|b(x)-b(y)||f(y)|\,dy\\
&\ge \frac{1}{r^{n-\alpha}}\int_{|y-x|\le r}|\Omega(x-y)||b(x)-b(y)||f(y)|\,dy.
\end{split}
\end{equation*}
Taking the supremum for all $r>0$ on both sides of the above inequality, we get
\begin{equation*}
[b,M_{\Omega,\alpha}](f)(x)\le[b,T_{|\Omega|,\alpha}](|f|)(x),\quad\mbox{for all}\; x\in\mathbb R^n.
\end{equation*}
Hence, as a direct consequence of Theorem 1.3, we finally obtain the following

\newtheorem{corollary}[theorem]{Corollary}
\begin{corollary}
Suppose that $\Omega\in L^s(S^{n-1})$ with $1<s\le\infty$ and $b\in BMO(\mathbb R^n)$. If $0<\alpha<n$, $1\le
s'<p<n/{\alpha}$, $1/q=1/p-{\alpha}/n$, $0<\kappa<p/q$ and $w^{s'}\in A(p/{s'},q/{s'})$, then the commutator
$[b,M_{\Omega,\alpha}]$ is bounded from $L^{p,\kappa}(w^p,w^q)$ to $L^{q,{\kappa q}/p}(w^q)$.
\end{corollary}

\end{document}